\documentclass{article}
\usepackage{latexsym,amsmath,amssymb,amscd}
\usepackage[dvips]{color}
\begin{document}

\title{On Rotating Axisymmetric Solutions of the Euler-Poisson Equations    }
\author{Juhi Jang \thanks{Department of Mathematics, University of Southern California, and Korea Institute for Advanced Study, Seoul, Republic of Korea, e-mail: juhijang@usc.edu} and Tetu Makino \thanks{Professor Emeritus at Yamaguchi University, Japan, e-mail: makino@yamaguchi-u.ac.jp}}
\date{\today}
\maketitle

\begin{abstract}
We consider stationary axisymmetric solutions of the Euler-Poisson equations, which govern the internal structure of barotropic gaseous stars. 
We take the general form of the equation of states which cover polytropic gaseous stars indexed by $6/5<\gamma <2$ and also white dwarfs. 
A generic condition of the existence of stationary solutions with differential rotation is given, and the existence of slowly rotating configurations near spherically symmetric equilibria is shown. The problem is formulated as a nonlinear integral equation, and is solved by an application of the infinite dimensional implicit function theorem. Oblateness of star surface is shown and also relationship between the central density and the total mass is given. \\


\end{abstract}

\newtheorem{Lemma}{Lemma}
\newtheorem{Proposition}{Proposition}
\newtheorem{Theorem}{Theorem}
\newtheorem{Definition}{Definition}
\newtheorem{Remark}{Remark}
\newtheorem{Corollary}{Corollary}

\numberwithin{equation}{section}

\section{Introduction}

We consider the Euler-Poisson equations
\begin{subequations}
\begin{align}
&\frac{\partial\rho}{\partial t}+\nabla\cdot (\rho v)=0, \label{1a} \\
&\rho\Big(\frac{\partial v}{\partial t}+
(v\cdot\nabla)v\Big)+\nabla P =-\rho\nabla\Phi, \label{1b}\\
&\triangle \Phi=4\pi\mathsf{G}\rho \label{1c}
\end{align}
\end{subequations}
on $(t,x)=(t, x_1, x_2, x_3) \in [0, T[ \times \mathbb{R}^3$. 
We assume \\

{\bf (A):}  {\it $P$ is a smooth function of $\rho >0$, $0<P, 0<dP/d\rho$ for
$\rho>0$, and there is a smooth function $\Lambda$ defined on $\mathbb{R}$ such that $\Lambda(0)=0$ and
\begin{equation}
P=\mathsf{A}\rho^{\gamma}(1+\Lambda(\mathsf{A}\rho^{\gamma-1}))
\quad\mbox{for}\quad \rho >0
\end{equation}
with positive constants $\mathsf{A}, \gamma$ such that $1<\gamma<2$}.\\

We put 
\begin{equation}
\nu:=\frac{1}{\gamma-1},
\end{equation}
while $1<\nu <+\infty$.

Supposing that the support of $\rho(t,\cdot)$ is compact, we replace \eqref{1c} by
\begin{equation}
\Phi(t,x)=-\mathsf{G}\int
\frac{\rho(t, x')}{|x-x'|}dx'. \label{4}
\end{equation}

Let us use the spherical polar co-ordinates
\begin{equation}
x_1=r\sqrt{1-\zeta^2}\cos\phi,\quad
x_2=r\sqrt{1-\zeta^2}\sin\phi,\quad
x_3=r\zeta. \label{5}
\end{equation}
We seek rotating axisymmetric stationary solutions
\begin{equation}
\rho=\rho(r,\zeta),\quad
v=-\Omega x_2\frac{\partial}{\partial x_1}+\Omega x_1
\frac{\partial}{\partial x_2}  \ \text{ \ on \ } \  \{\rho >0\} \label{1.6}
\end{equation}
where the angular velocity $\Omega$ is a given continuous function of $\varpi=r\sqrt{1-\zeta^2}=
\sqrt{(x_1)^2+(x_2)^2}$. 
Moreover we 
consider equatorially symmetric $\rho$, that is,
\begin{equation}
\rho(r, -\zeta)=\rho(r, \zeta)\quad\mbox{for}\quad \forall \zeta\in[-1,1].
\end{equation}
The support of the density is unknown a priori. This is a free boundary problem. When $\Omega=0$, the problem has spherically symmetric solutions  that satisfy the Lane-Emden equation. 

\

The study of self-gravitating rotating figures is a classical subject in celestial mechanics, astrophysics and mathematics. It goes back to Newton, Maclaurin, Jacobi, Poincar\'{e} and many others \cite{Chandrasekhar1967}: they studied the existence, oblateness, stability questions for uniformly rotating masses with the homogeneous density (incompressible). In the case of gas-like fluids, the compressibility must be taken into account, and much less is known even now.   

The first attempts to construct axisymmetric rotating solutions to the above compressible Euler-Poisson equations were made by astrophysicists E. A. Milne \cite{Milne}, H. von Zeipel \cite{Zeipel} in early 1920's and S. Chandrasekhar \cite{Chandra1933}, L. Lichtenstein \cite{Lichtenstein} in 1933 based on some perturbative methods.

A mathematically rigorous treatment of the problem was initiated by J. F. G. Auchmuty and R. Beals \cite{AuchmutyB} in 1971 using the formulation as variational problems to minimize the energy under the given total mass and the total angular momentum. Along this line, a lot of interesting progress has been made by many mathematicians on the existence, oblateness and stability \cite{Ach1991,CaF,ChLi,FrTur1,FrTur3,Li,LuoS0,LuoS}. 
However, such a variational method essentially requires $4/3<\gamma$. 

In \cite{JJTM}, we studied the problem by a non-variational approach, which is a natural justification of the perturbation method adopted by astrophysicists after S. Chandrasekhar \cite{Chandra1933}. We constructed slowly rotating solutions with small constant angular velocities  
 for $6/5<\gamma \leq 3/2$ by directly perturbing the unknown enthalpy function around the harmonic extension of spherically symmetric Lane-Emden solutions via the contraction mapping principle. 

After we completed our work  \cite{JJTM}, we were acquainted with the pioneering 
work by U. Heilig  \cite{Heilig} in 1994, which is a rediscovery of the old work by L. Lichtenstein \cite{Lichtenstein} in 1933. Their method is also non-variational, but the formulation is different from ours. They are seeking a mapping from the domain of the unperturbed configuration (e.g., the ball of spherically symmetric equilibrium) to the domain of perturbed configuration with the angular velocity which is near to that of the unperturbed one (e.g., small angular velocity) using the level surface of the density. Recently W. A. Strauss and Y.-L. Wu, \cite{StraussW} modified the framework of U. Heilig   \cite{Heilig} and  L. Lichtenstein \cite{Lichtenstein} to show the existence of rotating solutions with given total mass independent of the rotation speed.

The goal of this article is two-fold. The first one is to introduce a general implicit function theorem framework for our formulation where we directly perturb the enthalpy function as unknown. As a result, new solutions could be obtained not only near the Lane-Emden solutions but also near other nearby axisymmetric solutions. The second goal is to remove the restriction of $\gamma$ in the previous work \cite{JJTM} and to treat more general pressure laws. A new framework allows us to show the existence of slowly rotating stars for $6/5<\gamma<2$ including both uniform rotation and differential rotation, and also to show the existence of slowly rotating white dwarfs. 

As in many other applications of the implicit function theorem, a key is the invertibility of the main linear operator, which is equivalent to the triviality of the kernel. That was the key point in our previous work \cite{JJTM}, and we found that its parallelism played an important role also in the work by U. Heilig \cite{Heilig} and L. Lichtenstein \cite{Lichtenstein}. In this article, inspired by U. Heilig  and L. Lichtenstein, we give a new proof of the validity of the key kernel condition near the spherically symmetric equilibria. 

As a final remark of the introduction we would like to mention that our result yields a family of solutions of the same total mass. The central density is a free parameter that we fix in our construction. By a simple scaling argument, our solutions are parameterized by the central density $\rho_{\mathsf{O}}$ and the angular velocity squared $\Omega^2$. The total mass $M=\mathcal{M}(
\rho_{\mathsf{O}}, \Omega^2)$ is then determined for the unique configuration of $\rho_{\mathsf{O}}$ and $ \Omega^2$. In Section 6 we show that the central density is uniquely determined by the total mass for each given $\Omega^2$, provided that $\gamma \neq 4/3$ and $\Omega^2/\rho_{\mathsf{O}} $ is sufficiently small. In particular, it implies that for fixed total mass $\bar M = \mathcal{M}(\bar\rho_{\mathsf{O}},0)$, there exists a continuous curve $\mathcal C:\Omega^2 \rightarrow \rho_{\mathsf{O}}$ defined for sufficiently small $\Omega^2 \ll 1$ such that $ \mathcal{M}(\mathcal C(\Omega^2), \Omega^2)=\bar M$ and $\mathcal{C}(0)=\bar{\rho}_{\mathsf{O}}$.

This article proceeds as follows. Section 2 is devoted to the explanation of notations and the formulation of the problem. The problem is formulated as an integral equation as in \cite{JJTM}, and the Fr\'{e}chet derivative of the nonlinear integral operator to be considered is carefully analyzed under the wider condition $1<\nu< \infty$, where $\nu=1/(\gamma-1)$. (We discussed the case with $\nu \geq 2$ in the previous work \cite{JJTM}.) 
In Section 3 we shall state and prove the main results. A generic condition of the existence of differentially rotating configurations is given, and the existence of slowly rotating configurations near spherically symmetric equilibria is given.  In Section 4 we shall discuss a property of rotating configurations at the origin. In Section 5 we give the oblateness of the surface of slowly rotating configurations near spherically symmetric equilibria. In Section 6 we provide the relationship between the total mass and the central density, and prove that the central density is uniquely determined by the total mass as long as $\gamma\neq 4/3$. Lastly, in Section 7 we discuss the framework when the angular momentum is prescribed instead of the angular velocity.

\section{Notations and problem setting}

First of all we shall use the notations
$$ a\vee b :=\max(a, b), \qquad a\wedge b:=\min (a, b)$$
for numbers $a, b$ and we denote
$ u\vee 0 : \xi \mapsto u(\xi)\vee 0 $
for a function $u: \xi \mapsto u(\xi)$, that is,
$$
(u\vee 0 )(\xi)=
\begin{cases}
u(\xi) \quad \mbox{if}\quad u(\xi)\geq 0 \\
0\quad\mbox{if}\quad u(\xi) <0
\end{cases}.
$$\\

If a function $f$ of $r,\zeta$ such that $f(0,\zeta)=f(0,0) \forall \zeta \in [-1,1]$ is given, we define the function $f^{\flat}$ of $x=(x_1,x_2, x_3)$ 
by
$$f^{\flat}(x)=f(r,\zeta)$$
with \eqref{5}, that is, 
$$r=\sqrt{(x_1)^2+(x_2)^2+(x_3)^2},\qquad \zeta=\frac{x_3}{\sqrt{(x_1)^2+(x_2)^2+(x_3)^2}}.$$

We define 
\begin{equation}
\mathcal{K}f(r,\zeta)=\frac{1}{4\pi}
\int_{-1}^1\int_0^{\infty}
K(r,\zeta, r',\zeta')f(r',\zeta')r'^2dr'd\zeta' \label{2.1}
\end{equation}
with
\begin{equation}
K(r,\zeta,r',\zeta')=\int_0^{2\pi}
\frac{d\beta}{\sqrt{r^2+r'^2-2rr'(\sqrt{1-\zeta^2}\sqrt{1-\zeta'^2}\cos\beta+\zeta\zeta')}}
\end{equation}
while
\begin{equation}
\mathcal{K}^{\flat}g(x)=\frac{1}{4\pi}\int\frac{g(x')}{|x-x'|}dx'. \label{2.3}
\end{equation}
Then 
$$\mathcal{K}^{\flat}f^{\flat}=(\mathcal{K}f)^{\flat}$$
and
$$-\triangle(\mathcal{K}^{\flat}g)=g.$$

Given $R>0$, we shall denote
\begin{equation}
B(R)=\{ x\in \mathbb{R}^3 \  |\  |x|\leq R\}.
\end{equation}
The Euler-Poisson equations \eqref{1a},\eqref{1b},\eqref{4} under the ansatz \eqref{1.6} read
\begin{subequations}
\begin{align}
-\rho(1-\zeta^2)r\Omega^2+\frac{\partial P}{\partial r}&=-\rho
\frac{\partial\Phi}{\partial r}, \label{8a} \\
\rho\zeta r^2\Omega^2+\frac{\partial P}{\partial\zeta}&=
-\rho\frac{\partial \Phi}{\partial\zeta} \label{8b} \\
\Phi&=-4\pi\mathsf{G}\mathcal{K}\rho.
\end{align}
\end{subequations}

We introduce the enthalpy variable $u$ by
\begin{equation}
u=\int_0^{\rho}\frac{dP}{\rho}.
\end{equation}

Then
\begin{align}
&u=\frac{\gamma\mathsf{A}}{\gamma-1}\rho^{\gamma-1}
(1+\Lambda_u(\mathsf{A}\rho^{\gamma-1})), \\
\Leftrightarrow &  \nonumber \\
&\rho=\Big(\frac{\gamma-1}{\gamma\mathsf{A}}\Big)^{\nu}
u^{\nu}(1+\Lambda_{\rho}(u)).
\end{align}
Here
$$\Lambda_u(\xi)=\frac{1}{\xi}\int_0^{\xi}
\Lambda_1(\xi')d\xi',
$$
while
$$\Lambda_1(\xi)=\Big[1+\frac{\gamma-1}{\gamma}\xi
\frac{d}{d\xi}\Big]\Lambda(\xi).
$$
We see $\Lambda_u(0)=\Lambda_1(0)=0$.

The system \eqref{8a},\eqref{8b} reduces to
\begin{subequations}
\begin{align}
-(1-\zeta^2)r\Omega^2+\frac{\partial u}{\partial r}&=
-\frac{\partial\Phi}{\partial r}, \\
\zeta r^2\Omega^2+\frac{\partial u}{\partial\zeta}&=
-\frac{\partial\Phi}{\partial\zeta},
\end{align}
\end{subequations}
on $\{\rho >0\}$, which is equivalent to
\begin{equation}
\Phi+u=
\mathsf{B}(r\sqrt{1-\zeta^2})+\mbox{Const.},
\end{equation}
where
\begin{equation}
\mathsf{B}({\varpi}):=\int_0^{\varpi}\Omega(\varpi')^2\varpi'd\varpi'.
\end{equation}

Let us introduce the parameters $u_{\mathsf{O}}>0$, which should mean 
the central enthalpy $u(0,0)$, and 
\begin{equation}
\mathsf{a}=\frac{1}{\sqrt{4\pi\mathsf{G}}}
\Big(\frac{\mathsf{A}\gamma}{\gamma-1}\Big)^{\frac{1}{2(\gamma-1)}}
u_{\mathsf{O}}^{-\frac{2-\gamma}{2(\gamma-1)}}, 
\end{equation}
 and let us put
\begin{equation}
\mathsf{b}(\varpi)
=
u_{\mathsf{O}}^{-1}\mathsf{B}(\mathsf{a}\varpi).
\end{equation}
Then by the scale change $x \mapsto x/\mathsf{a}, r \mapsto r/\mathsf{a}, u\mapsto 
u/u_{\mathsf{O}}$, the problem reduces to the integral equation
\begin{equation}
u=\mathfrak{g}+\mathcal{G}(u), \label{E}
\end{equation}
where
\begin{align}
&\mathfrak{g}(r,\zeta)=\mathsf{b}(r\sqrt{1-\zeta^2}) \\
&\mathcal{G}(u)=1+\mathcal{K}(\mathsf{f}\circ u)-
\mathcal{K}(\mathsf{f}\circ u)(0,0), \\
&\mathsf{f}(u)=(u\vee 0)^{\nu}(1+\Lambda_{\rho}(u_{\mathsf{O}}u)). \label{Def.f}
\end{align}\\

Let $r_{\infty}>0$. 

Since 
$$\nabla \mathfrak{g}^{\flat}=
\frac{x_1}{\varpi}\frac{d\mathsf{b}}{d\varpi}\frac{\partial}{\partial x_1}
+
\frac{x_2}{\varpi}\frac{d\mathsf{b}}{d\varpi}\frac{\partial}{\partial x_2},
$$ we see that $\nabla\mathfrak{g}^{\flat} \in C(B(r_{\infty}))$ if
$\mathsf{b} \in C^1([0,r_{\infty}])$ and $\displaystyle \frac{d\mathsf{b}}{d\varpi}=0$ at $\varpi =0$. So, we consider the Banach space
\begin{equation}
\mathfrak{B}:=\{ \mathsf{b}\in C^1([0,r_{\infty}])\  |\  
\frac{d\mathsf{b}}{d\varpi}=0 \quad\mbox{at}\quad \varpi=0 \}
\end{equation}
endowed with the norm
\begin{equation}
\|\mathsf{b}\|_{\mathfrak{B}}=\|\mathsf{b}\|_{L^{\infty}(0,r_{\infty})}+
\Big\|\frac{d\mathsf{b}}{d\varpi}\Big\|_{L^{\infty}(0,r_{\infty})}.
\end{equation}
Note that if $\Omega \in C([0,\mathsf{a}r_{\infty} ])$ then $\mathsf{b} \in \mathfrak{B}$ and
\begin{equation}
\|\mathsf{b}\|_{\mathfrak{B}}\leq \frac{1}{4\pi\mathsf{G}}
\Big(\frac{\mathsf{A}\gamma}{\gamma-1}\Big)^{\frac{1}{\gamma-1}}
u_{\mathsf{O}}^{-\frac{1}{\gamma-1}}
\Big(\frac{r_{\infty}^2}{2}+r_{\infty}\Big)\|\Omega\|_{L^{\infty}(0, \mathsf{a}r_{\infty} )}^2.
\end{equation}

We consider the integral equation \eqref{E} in the Banach space
\begin{align}
\mathfrak{E}&=\{ u\in C([0,r_{\infty}]\times [-1,1] | u(0,\zeta)=u(0,0), \nonumber \\
&u(r,-\zeta)=u(r,\zeta) \forall \zeta \in [-1,1] \forall
r\in [0,r_{\infty}]\}
\end{align}
endowed with the norm
\begin{equation}
\|u\|_{\mathfrak{E}}=\|u\|_{L^{\infty}([0,r_{\infty}]\times [-1,1])}
\end{equation}

Given $u \in \mathfrak{E}$, for which the support of $u^{\flat}$ may not be a compact subset of the interior of the ball $B(r_{\infty})$, we define $\mathcal{G}(u)$ by
\begin{equation}
\mathcal{K}(\mathsf{f}\circ u)(r,\zeta)=\frac{1}{4\pi}
\int_{-1}^1\int_0^{r_{\infty}}
K(r,\zeta, r',\zeta')\mathsf{f}(u(r',\zeta'))r'^2dr'd\zeta' 
\end{equation}
instead of \eqref{2.1} with $f=\mathsf{f}\circ u$.

Then the Fr\'{e}chet derivative $D\mathcal{G}(u)$ is given by
\begin{equation}
D\mathcal{G}(u)h=\mathcal{K}(\mathsf{f}'\circ u\cdot h)
-\mathcal{K}(\mathsf{f}'\circ u \cdot h)(0,0),
\end{equation}
where
\begin{align}
\mathsf{f}'(u)&=\nu (u\vee 0)^{\nu-1}
(1+\Lambda_{\rho'}(u_{\mathsf{O}}u)), \label{2.25}\\
\Lambda_{\rho'}(\xi)&=\Big[1+\frac{1}{\nu}\xi\frac{d}{d\xi}\Big]\Lambda_{\rho}(\xi).\label{2.26}
\end{align}
We see $\Lambda_{\rho'}(0)=0$.

  See \cite[Proposition 1]{JJTM}. Actually \cite[Lemma 3]{JJTM} is valid even if $\nu <2$, that is,

   \begin{Lemma} \label{Lem.1}
    Let $1<\nu<\infty$.
Let $ M<\infty$ be fixed and let $|u|, |u+h|\leq M$, and put
\begin{equation}
((u+h)\vee 0)^{\nu}-(u \vee 0)^{\nu}=
\nu (u \vee 0)^{\nu-1}h+\mathfrak{R}(u,h).
\end{equation}
Then there is a constant C depending on $M$ such that
\begin{equation}
|\mathfrak{R}(u,h)|\leq C|h|^{\nu\wedge 2}.
\end{equation}
\end{Lemma}

Proof. It is easy to verify
$$|(1+x)^{\nu-1}-1|\lesssim
\begin{cases}
|x|\quad\mbox{for}\quad |x|\leq 1 \\
x^{\nu-1}\quad\mbox{for}\quad 1\leq x.
\end{cases}
$$
Then it follows that
$$|(1+x)^{\nu}
-1-
\nu x|\lesssim
\begin{cases}
|x|^2\quad\mbox{for}\quad |x|\leq 1 \\
x^{\nu}\quad\mbox{for}\quad 1\leq x.
\end{cases}
$$

 Now we consider
$$\Delta\rho:=((u+h)\vee 0)^{\nu}-(u \vee 0)^{\nu}.
$$

Case-(00): Suppose $u>0$ and $u+h>0$. Then
\begin{align*}
\Delta\rho&=(u+h)^{\nu}-u^{\nu}=u^{\nu}\Big(\Big(1+
\frac{h}{u}\Big)^{\nu}-1\Big)\Big)= \\
&=\begin{cases}
u^{\nu}\Big(
\nu \frac{h}{u}+O\Big(\Big|\frac{h}{u}\Big|^2\Big)\Big) 
\quad\mbox{for}\quad |h|\leq u \\
u^{\nu}\Big(
\nu\frac{h}{u}+O\Big(\Big|\frac{h}{u}\Big|^{\nu}\Big)\Big)
\quad\mbox{for}\quad u\leq h
\end{cases}\\
&=\begin{cases}
\nu u^{\nu-1}h+O(u^{\nu-2}|h|^2)\quad\mbox{for}\quad |h|\leq u \\
\nu u^{\nu-1}h+O(|h|^{\nu})\quad\mbox{for}\quad u\leq h
\end{cases}\\
&=\nu u^{\nu-1}h+O(|h|^{\nu\wedge 2}),
\end{align*}
where we have used the following estimates:
Let $|h| \leq u$; When $2\leq\nu$, then
$$u^{\nu-2}|h|^2\leq M^{\nu-2}|h|^2; $$
When $1<\nu <2$, then
$$u^{\nu-2}|h|^2 =\Big|\frac{h}{u}\Big|^{2-\nu}|h|^{\nu}
\leq |h|^{\nu}.$$

Case-(01): Suppose $u>0$ but $u+h\leq 0$. Then
$\Delta\rho=-u^{\nu}$. But $0<u\leq -h$ implies
$|u|\leq |h|$ and
$$|\Delta\rho|\leq |h|^{\nu}\lesssim |h|^{\nu\wedge 2},$$
provided that $|h|\leq|u+h|+|u|\leq2 M$ when $\nu > 2$ and $\nu\wedge 2=2$. On the other hand,
$$|\nu (u \vee 0)^{\nu-1}h|\leq \nu u^{\nu-1}|h|\leq\nu |h|^{\nu}.$$

Case-(10): Suppose $u\leq 0$ but $u+h>0$. Then
$\Delta\rho=(u+h)^{\nu}, 0<u+h\leq h$,
and $(u \vee 0)=0$.

Case-(11): Suppose $u\leq 0$ and $u+h\leq 0$. Then
$\Delta\rho=0$ and $(u \vee 0)=0$, no problem. $\square$ \\

  Similarly we have
 
\begin{Lemma}\label{Lem.2} 
Let $1<\nu<\infty$. Let $M<\infty$ be fixed and let $|u|,|u+h| \leq M$. Then
\begin{equation}
|((u+h)\vee 0)^{\nu-1}-(u\vee 0)^{\nu-1}|\leq C|h|^{(\nu-1)\wedge 1}.
\end{equation}
Here $C$ is a constant depending on $M$.
\end{Lemma}

Proof. Put
$$\Delta q=((u+h)\vee 0)^{\nu-1}-(u\vee 0)^{\nu-1}.$$
Case-(00): Suppose $u>0, u+h>0$. Then
\begin{align*}
\Delta q &=u^{\nu-1}\Big(\Big(1+\frac{h}{u}\Big)^{\nu-1}-1\Big) \\
&\lesssim u^{\nu-1}\times
\begin{cases}
\Big|\frac{h}{u}\Big| \quad\mbox{for}\quad \Big|\frac{h}{u}\Big|\leq 1 \\
\Big(\frac{h}{u}\Big)^{\nu-1}\quad\mbox{for}\quad \frac{h}{u}\geq 1
\end{cases}\\
&=\begin{cases}
u^{\nu-2}|h|\quad\mbox{for}\quad |h|\leq u \\
h^{\nu-1}\quad\mbox{for}\quad h\geq u.
\end{cases}
\end{align*}
But, when $1<\nu \leq 2$, then 
$$u^{\nu-2}|h|=\Big(\frac{1}{u}\Big)^{2-\nu}|h| \leq \Big(\frac{1}{|h|}\Big)^{2-\nu}|h|=|h|^{\nu-1}, $$
provided that $|h|\leq u$; when $2\leq \nu$, then
$$u^{\nu-2}|h|\leq M^{\nu-2}|h|=M^{\nu-2}|h|^{(\nu-1)\wedge1}.$$

Case-(01): Suppose $u>0, u+h\leq 0$. Then $\Delta q=-u^{\nu-1}$. But
$0<u\leq -h$ implies $|u|\leq |h|$ and
$|\Delta q|\leq |h|^{\nu-1}$. 

Case-(10): Suppose $u\leq 0, u+h>0$. Then
$\Delta q=(u+h)^{\nu-1}, 0<u+h\leq h$ sot hat
$0< \Delta q \leq h^{\nu-1}$.

Case-(11): Suppose $u \leq 0, u+h\leq 0$. Then $\Delta q=0$, no problem. $\square$

  \section{Results}
  
  Let us fix $  r_{\infty} > 0$ . We suppose that 
a function  $\mathsf{b} \in \mathfrak{B}$ and a positive number $u_{\mathsf{O}}$ are given. ( In order to give $\mathsf{b} \in \mathfrak{B}$, it is sufficient to give  $\Omega \in C([0,r_{\infty}])$.)
  
  We are considering the integral equation \eqref{E} on the Banach space $\mathfrak{E}$. We shall find a solution $u \in \mathfrak{E}$ which satisfies the following two conditions:\\
   
  {\bf (a1):} 
  {\it It holds
$$\frac{\partial u}{\partial r} <0$$
for $r_0\leq r \leq r_{\infty}, |\zeta|\leq 1$,
$r_0$ being a small positive number;}\\
  
{\bf (a2): }   
  {\it There exists a continuous function
  $R(\zeta)$ of $|\zeta|\leq 1$ such that
  $r_0<R(\zeta) <r_{\infty} \quad \forall \zeta \in [-1,1]$ and
  \begin{align*}
  &u(r,\zeta) >0\quad\mbox{for}\quad 0\leq r<R(\zeta) \\
  &u(r,\zeta)=0\quad\mbox{for}\quad r=R(\zeta) \\
 & u(r,\zeta)<0\quad\mbox{for}\quad R(\zeta)<r\leq r_{\infty};
  \end{align*}
Here $r_0$ is the positive number in} {\bf (a1)}\\

Here we have noted that the continuous derivative $\partial u/\partial r$ exists on
$0<r\leq r_{\infty}, |\zeta|\leq 1$. Actually $u^{\flat}$ satisfies 
$$u^{\flat}=\mathfrak{g}^{\flat}
+1+\mathcal{K}^{\flat}(\mathsf{f}\circ u^{\flat})
-\mathcal{K}^{\flat}(\mathsf{f}\circ u^{\flat})(O),$$
Here we consider
$$\mathcal{K}^{\flat}(\mathsf{f}\circ u^{\flat})(x)
=\frac{1}{4\pi}\int_{|x'|\leq r_{\infty}}
\frac{\mathsf{f}(u^{\flat}(x'))}{|x-x'|}dx' $$
instead of \eqref{2.3} with $g=\mathsf{f}\circ u^{\flat}$.
 Therefore as well-known
$\mathcal{K}^{\flat}(\mathsf{f}\circ u^{\flat}) \in C^1(\mathbb{R}^3)$ and
$$|\nabla\mathcal{K}^{\flat}(\mathsf{f}\circ u^{\flat})|\leq C\|u\|_{\mathfrak{E}}^{\nu}.$$
Since $\mathfrak{g}^{\flat}\in C^1$, this implies $u^{\flat}\in C^1(B(r_{\infty}))$
and
$$\frac{\partial u}{\partial r}=\frac{x\cdot\nabla u^{\flat}}{r}
\in C(0<r\leq r_{\infty}, |\zeta|\leq 1).$$

\begin{Definition}
We shall call a solution $u\in\mathfrak{E}$ of \eqref{E} which satisfies the above conditions {\bf (a1) (a2)} an admissible solution of \eqref{E} with parameters 
$\mathsf{b}, u_{\mathsf{O}}.$
\end{Definition}

Let $u$ be an admissible solution. Put
$$\mathfrak{R}=\{ x \in B(r_{\infty}) | \ u^{\flat} (x)>0\}
$$
and
$$S=\partial\mathfrak{R}=\{ x | u^{\flat}(x)=0\}.$$
Then $S$ is a compact $C^1$-surface in the Euclidean space
$\mathbb{R}^3$, since $u^{\flat}\in C^1$. 

The physical vacuum boundary condition is
$$-\infty <\frac{\partial u^{\flat}}{\partial N}\Big|_S <0,$$
where $N$ is the unit normal vector outer with respect to the domain $\mathfrak{R}=\{u^{\flat}>0\}$. So,
$$\frac{\partial u^{\flat}}{\partial N}=-|\nabla u^{\flat}|. $$
But we see
$$|\nabla u^{\flat}|^2=\Big(\frac{\partial u}{\partial r}\Big)^2+
\frac{1-\zeta^2}{r^2}\Big(\frac{\partial u^{\flat}}{\partial\zeta}\Big)^2 >0$$
along $S$, since 
$\partial u/\partial r <0$ on $r_0\leq r\leq r_{\infty}$, and
$r_0 < R(\zeta) <r_{\infty}$. That is, any admissible solution satisfies the physical vacuum boundary condition. 

We consider the following Heilig-Lichtenstein condition {\bf (HL)}:

\begin{Definition}\label{Def}
Let $r_{\infty}>0$ be fixed. A function $u \in \mathfrak{E}$ is said to satisfy the  condition {\bf (HL)} if
\begin{equation}
h \in \mathfrak{E}\quad\mbox{and}\quad h=D\mathcal{G}({u})h
\qquad\Rightarrow \qquad  h=0. 
\end{equation}
\end{Definition}

We claim the following
\begin{Theorem}\label{Th.1}
Let $r_{\infty}>0$. Let $\bar{u} \in\mathfrak{E}$ be an admissible solution of \eqref{E} with parameters $\bar{\mathsf{b}}, \bar{u}_{\mathsf{O}}$. Here $\bar{\mathsf{b}} \in \mathfrak{B}$ and
$\bar{u}_{\mathsf{O}} >0$.
Suppose that $u=\bar{u}$ enjoys the condition {\bf (HL) }. Then there exist small positive constants $\delta, \epsilon$ such that
for $(\mathsf{b} , u_{\mathsf{O}} ) \in \mathfrak{B}\times \mathbb{R} $ such that 
$\|\mathsf{b}-\bar{\mathsf{b}}\|_{\mathfrak{B}}
+|u_{\mathsf{O}}-\bar{u}_{\mathsf{O}}| \leq\delta$
there is a unique admissible solution $u$ with parameters 
$\mathsf{b}, u_{\mathsf{O}}$ which enjoys {\bf (HL)} such that $\| u-\bar{u}\|_{\mathfrak{E}}\leq\epsilon$, and $(\mathsf{b}, u_{\mathsf{O}}) \mapsto u$ is continuous.
\end{Theorem}

In order to prove Theorem \ref{Th.1} we shall apply the following implicit function theorem: \\

{\bf Implicit function theorem } (\cite[Theorem 15.1]{Deimling}):
{\it Let $X,Y,Z$ be Banach spaces, $U \subset X$ and $V\subset Y$ neighborhood of $x_0$ and $y_0$ respectively., $ F:U\times
V \rightarrow Z$ continuous and continuously differential with respect to $y$. Suppose also that $F(x_0, y_0)=0$ and
$F_y^{-1}(x_0,y_0) \in L(Z,Y)$. Then there exists balls $\bar{B}_r(x_0)\subset U,
\bar{B}_{\delta}(y_0)\subset V$ and exactly one map
$T:B_r(x_0) \rightarrow B_{\delta}(y_0)$ such that 
$Tx_0=y_0$ and $F(x, Tx)=0$ on $B_r(x_0)$. This map $T$ is continuous.
}\\

Proof of Theorem \ref{Th.1}. The equation \eqref{E} can be written as 
\begin{equation}
\mathfrak{F}(\mathsf{b}, u_{\mathsf{O}}, u)=0,\label{3.1}
\end{equation}
where
\begin{equation}
\mathfrak{F}(\mathsf{b}, u_{\mathsf{O}}, u)=u-\mathfrak{g}-\mathcal{G}(u).
\end{equation}
Then  $\mathfrak{F}(\mathsf{b},u_{\mathsf{O}}, u)$ and
$D_u\mathfrak{F}(\mathsf{b},u_{\mathsf{O}}, u)$ are continuous on
$$
\{ (\mathsf{b}, u_{\mathsf{O}},u) \in \mathfrak{B}\times \mathbb{R} 
\times \mathfrak{E}\  | \  
\|\mathsf{b}-\bar{\mathsf{b}}\|_{\mathfrak{B}}+
|u_{\mathsf{O}}-\bar{u}_{\mathsf{O}}|+\|u-\bar{u}\| \leq \delta_0 \ll 1\}.$$
Here
\begin{align}
&D_u\mathfrak{F}(\mathsf{b},u_{\mathsf{O}}, u)=I-D\mathcal{G}(u), \\
&D\mathcal{G}(u)h=
\mathcal{K}(\mathsf{f}'\circ u\cdot h)-
\mathcal{K}(\mathsf{f}'\circ u\cdot h)(0,0).
\end{align}

As for the continuity of $D_u\mathfrak{F}$ with respect to $u$, see Lemma \ref{Lem.2}.

Hence if the condition {\bf (HL)} holds, then
$\mbox{Ker}(I-D\mathcal{G}(\bar{u}))=\{0\}$. Since $D\mathcal{G}(\bar{u})$ is a compact operator in 
$\mathfrak{E}$, this implies that there exists the bounded linear inverse 
$(I-D\mathcal{G}(\bar{u}))^{-1}$. Then 
the implicit function theorem \cite[Theorem 15.1]{Deimling} can be applied:
for $\|\mathsf{b}-\bar{\mathsf{b}}\|_{\mathfrak{B}}+|u_{\mathsf{O}}-\bar{u}_{\mathsf{O}}|\leq \delta$ there exists a unique solution $u
\in \mathfrak{E}$ of \eqref{3.1} such that $\|u-\bar{u}\|_{\mathfrak{E}}\leq\epsilon$. Then we have
$$|\nabla u^{\flat}-\nabla\bar{u}^{\flat}|
\leq |\nabla \mathfrak{g}^{\flat}-\nabla\bar{\mathfrak{g}}^{\flat}|+
C\|u-\bar{u}\|_{\mathfrak{E}} \leq C'(\delta+\epsilon).$$
Since we can assume that $\delta'=C'(\delta+\epsilon)$ is sufficiently small, we can assume that 
$$\Big|\frac{\partial u}{\partial r}-\frac{\partial\bar{u}}{\partial r}\Big|
\leq |\nabla u^{\flat}-\nabla \bar{u}^{\flat}|\leq \delta'$$
implies the condition {\bf (a1)}, since
$$\frac{\partial\bar{u}}{\partial r} \leq -\frac{1}{C}<0\quad\mbox{uniformly for}\quad r_0\leq r\leq r_{\infty}, |\zeta|\leq 1.$$
Moreover the implicit function theorem guarantees that for any 
$\zeta \in [-1,1]$ the equation $u(r,\zeta)=0$ admits a unique continuous solution $r=R(\zeta)$ such that 
$$|R(\zeta)-\bar{R}(\zeta)|\leq \epsilon', 
\quad r_0<R(\zeta)<r_{\infty}, $$
$\epsilon'$ being small. Thus the condition {\bf (a2)} holds.

Let us verify that $u$ enjoys the condition {\bf (HL)}, provided that $\epsilon$ is sufficiently small. The bounded inverse
$(I-D\mathcal{G}(\bar{u}))^{-1}$ exists and the solution $u$ satisfies
\begin{align*}
\|(D\mathcal{G}(u)-\mathcal{G}(\bar{u}))h\|_{\mathfrak{E}} &
=\|\mathcal{K}(\mathsf{f}'\circ u -\mathsf{f}'\circ \bar{u})\cdot h\|_{\mathfrak{E}} \\
&\leq C\|u-\bar{u}\|_{\mathfrak{E}}^{ (\nu-1)\wedge 1 }\|h\|_{\mathfrak{E}},
\end{align*}
thanks to Lemma \ref{Lem.2}. Therefore the operator
$$\mathcal{U}:=(I-D\mathcal{G}(\bar{u}))^{-1}
(D\mathcal{G}(u)-D\mathcal{G}(\bar{u}))$$
enjoys the estimate of the operator norm 
$$|\|\mathcal{U}\||_{\mathfrak{E};\mathfrak{E}} \leq \epsilon' <1,$$
provided that $\|u-\bar{u}\|_{\mathfrak{E}}\leq\epsilon$ is sufficiently small. Then the Neumann series gives the existence of the bounded inverse
$$(I-D\mathcal{G}(u))^{-1}=
(I-\mathcal{U})^{-1}
(I-D\mathcal{G}(\bar{u}))^{-1}, $$
a fortiori, the condition {\bf (HL)} for $u$.

This completes the proof. $\square$

\begin{Theorem}\label{Th.2}
Let $u_{\mathsf{O}} >0$ be fixed. Suppose that the solution  $u(r)$ of the initial value problem of the ordinary differential equation
\begin{align}
&-\frac{1}{r^2}\frac{d}{dr}r^2\frac{du}{dr}=\mathsf{f}(u), \label{SE}\\
&u=1+O(r^2)\quad\mbox{as}\quad r\rightarrow +0
\end{align}
exists to be positive on $0<r<R$ and $u(r)\rightarrow 0$ as $r\rightarrow R-0$. Choose $ r_{\infty}$ such that
$R<r_{\infty}$. Then we can extend $u(r)$ onto $[0, r_{\infty}]$
so that the extension is an admissible solution with parameters  $\mathsf{b}=0, u_{\mathsf{O}}$ which satisfies the condition {\bf (HL)}.
\end{Theorem}

Proof. We define $u(r)$ for $r\geq R$ by
$$u(r)=-\mu\Big(\frac{1}{R}-\frac{1}{r}\Big),$$
where
$$\mu=\int_0^R\mathsf{f}(u(r))r^2dr=-r^2\frac{du}{dr}\Big|_{r=R-0}>0.$$
Then this extension $u(r)$ satisfies \eqref{SE}, $du/dr <0$ for $0<r<+\infty$ and $u(r)=O(1/r), du/dr =O(1/r^2)$ as $r \rightarrow +\infty$. Thus we have an admissible solution. 

We are going to verify the condition {\bf (HL)}. Let $h\in\mathfrak{E}$ satisfy
\begin{align*}
h=&D\mathcal{G}(u)h \\
=&\mathcal{K}(\mathsf{f}'\circ u\cdot h)
-\mathcal{K}(\mathsf{f}'\circ u\cdot h)(0,0).
\end{align*}
This implies 
$$-\triangle h=\mathsf{f}'(u(r))h.$$
We can assume that $h(r,\zeta)$ exists and satisfies this equation on $0 \leq 
r<+\infty, |\zeta|\leq 1$ and $h=O(1/r), \partial h/\partial r =O(1/r^2)$ as $r\rightarrow +\infty$,
since $\mathsf{f}'(u(r))$ vanishes on $r >R$. 

Let us consider 
$$h_j(r)=\frac{2j+1}{2}\int_{-1}^1h(r,\zeta)P_j(\zeta)d\zeta,$$
where $P_j(\xi)$ is the Legendre's polynomial. We want to show $h_j(r)$ vanishes for any even $j$. Note that  $y=h_j(r)$ is the solution of 
$$\Big[-\frac{1}{r^2}\frac{d}{dr}r^2\frac{d}{dr}+
\frac{j(j+1)}{r^2}\Big]y=
\mathsf{f}'(u(r))y, \eqno{(E_j)}$$
$$y|_{r=0}=0, \quad \frac{dy}{dr}\Big|_{r=0} =O(1). $$

When $j=0$, $(E_0)$ has a fundamental system of solutions 
$y_1 \sim 1,  y_2\sim  \displaystyle \frac{1}{r}$ as $r \rightarrow +0$. So, we see $h_0=0$. 

Let $j\geq 2$. We use the trick of \cite[Drittes Kapitel, \S 3]{Lichtenstein},
\cite[Section 5]{Heilig}. Since $h_j(r)=O(1/r), dh_j(r)/dr =O(1/r^2)$ as $r\rightarrow +\infty$, we have the identity
\begin{align}
h_j(r)&=\frac{1}{2j+1}\Big[
\frac{1}{r^2}\int_0^r
\mathsf{f}'(u(r'))h_j(r')\Big(\frac{r'}{r}\Big)^{j-1}r'^3dr'+ \nonumber \\
&+r\int_r^R
\mathsf{f}'(u(r'))h_j(r')\Big(\frac{r}{r'}\Big)^{j-1}dr'\Big]
\end{align}
for $0<r \leq R$.  See Appendix, namely, note that we apply Proposition \ref{Prop.A} to $A=0$. 

Put
\begin{equation}
\psi(r)=\frac{1}{r^2}\int_0^r\mathsf{f}(u(r'))r'^2dr'=-\frac{du}{dr},
\end{equation}
and consider 
\begin{equation}
H(r)=\frac{h_j(r)}{\psi(r)}.
\end{equation}
Note that
$$\mathsf{f}'(u(r))\cdot\psi(r)=-\frac{d\mathsf{f}(u(r))}{dr} >0,$$
since $\mathsf{f}'(u) \propto d\rho/du =\rho d\rho/dP >0$.
Hence
\begin{align*}
H(r)&=\frac{1}{(2j+1)\psi(r)}\Big[
\frac{1}{r^2}\int_0^r
\mathsf{f}'(u(r'))\psi(u(r'))H(r')\Big(\frac{r'}{r}\Big)^{j-1}r'^3dr' + \\
&+r\int_r^R
\mathsf{f}'(u(r'))\psi(r')H(r')\Big(\frac{r}{r'}\Big)^{j-1}dr'\Big]
\end{align*}
implies
\begin{align*}
|H(r)|&\leq\frac{1}{(2j+1)\psi(r)}\Big[
\frac{1}{r^2}\int_0^r-\frac{d\mathsf{f}(u(r'))}{dr'}r'^3dr' 
+r\int_r^R-\frac{d\mathsf{f}(u(r'))}{dr'}dr'\Big]\|H\|_{L^{\infty}} \\
&=
\frac{1}{(2j+1)\psi(r)}\Big[
\frac{1}{r^2}\int_0^r-\frac{d\mathsf{f}(u(r'))}{dr'}r'^3dr' 
+r
\mathsf{f}(u(r))\Big]\|H\|_{L^{\infty}} \\
&=\frac{3}{(2j+1)\psi(r)}\Big[
\frac{1}{r^2}\int_0^r\mathsf{f}(u(r'))r'^2dr'\Big]\|H\|_{L^{\infty}} \\
&=\frac{3}{2j+1}\|H\|_{L^{\infty}},
\end{align*}
that is, 
$$\|H\|_{L^{\infty}}\leq \frac{3}{2j+1}\|H\|_{L^{\infty}}.$$
Note that
 $H$ is bounded, since $|h_j(r)|\leq Cr$ and
$\psi(r)\geq r/C, C$ being a positive constant. But
$\displaystyle \frac{3}{2j+1}<1$ for $j\geq 2$, so ,$ \|H\|_{L^{\infty}}=0$,
that is, $h_j(r)=0$. $\square$.\\

Combining Theorem \ref{Th.1} and Theorem \ref{Th.2}, we have

\begin{Corollary}\label{Cor.1}
Let $\bar{u}_{\mathsf{O}} >0$ be fixed. Suppose that the solution  $\bar{u}(r)$ of the initial value problem of the ordinary differential equation
\begin{align}
&-\frac{1}{r^2}\frac{d}{dr}r^2\frac{d\bar{u}}{dr}=\bar{\mathsf{f}}(\bar{u}), \\
&\bar{u}=1+O(r^2)\quad\mbox{as}\quad r\rightarrow +0
\end{align}
exists to be positive on $0<r<\bar{R}$ and $u(r)\rightarrow 0$ as $r\rightarrow \bar{R}-0$. 
Here 
$$\bar{\mathsf{f}}(u)=(u\vee 0)^{\nu}(1+\Lambda_{\rho}(\bar{u}_{\mathsf{O}}u))).$$
Choose $ r_{\infty}$ such that
$\bar{R}<r_{\infty}$. Then there are positive positive numbers $\delta $ and $\epsilon$ such that for $(\mathsf{b}, u_{\mathsf{O}}) \in \mathfrak{B}
\times \mathbb{R}$ with 
$\|\mathsf{b}\|_{\mathfrak{B}}+|u_{\mathsf{O}}-\bar{u}_{\mathsf{O}}|\leq \delta$ there is a unique admissible solution $u$ with parameters $\mathsf{b}, u_{\mathsf{O}}$ which enjoys {\bf (HL)} such that $\|u-\bar{u}\|_{\mathfrak{E}}\leq \epsilon$,
where $\bar{u}$ is extended to be harmonic on $r \geq \bar{R}$, and $(\mathsf{b}, u_{\mathsf{O}}) \mapsto u$ is continuous. 

\end{Corollary}

{\bf Applications:}\\

\textbullet \  Particularly, when $P=\mathsf{A}\rho^{\gamma}, 6/5<\gamma<2$,
then the Lane-Emden equation
$$-\frac{1}{r^2}\frac{d}{dr}r^2\frac{du}{dr}=u^{\nu},\quad u=1+O(r^2)
\quad(r\rightarrow +0) $$
admits the Lane-Emden function $\theta(r;\nu)$ with finite zero $\xi_1(\nu)$. 
See \cite{Chandra} and \cite{JosephL}. Take $r_{\infty}>\xi_1(\nu)$, Theorem \ref{Th.1}, Theorem \ref{Th.2}, or Corollary \ref{Cor.1},  can be applied, and we extend the result of \cite{JJTM} obtained for $6/5 <\gamma \leq 3/2$. \\

\textbullet \  Also, the existence of slowly rotating white dwarfs is given. Actually the equation of state for white dwarfs reads
$$P=Ac^5F(X),\qquad \rho=Bc^3X^3, $$
with
$$F(X)=X(2X^3-3)\sqrt{X^2+1}+3\mbox{sinh}^{-1}X,$$
$A,B, c$ are positive constants, so that $\mathsf{A}=8A/5B^{5/3}, \gamma=5/3$ and
$$
\rho=\frac{B^{5/2}}{4A^{3/2}}u^{\frac{3}{2}}\Big(1+
\frac{B}{16A}\frac{u}{c^2}\Big)^{\frac{3}{2}}.$$
See \cite[Chapter XI]{Chandra}.  Since $\gamma >4/3$ ( or $\nu <3$ ), we can claim that for any central density the radius
of the spherically symmetric white dwarf  is finite.
 For a rigorous proof of this claim, see \cite{TM84}.
 
 \section{Monotonicity at the origin}
 
 Recall that for an admissible solution $u$ we require that {\bf (a1)}: 
 $$\frac{\partial u}{\partial r} <0\quad\mbox{for}\quad
 r_0\leq r\leq r_{\infty}, |\zeta|\leq 1,$$
 $r_0$ being a number such that $0<r_0<r_{\infty}$. However if $\mathsf{b}=0$ and $u$ is
 spherically symmetric, then we have
 $$\frac{du}{dr}=-\frac{1}{r^2}\int_0^r
 \mathsf{f}(u(r'))r'^2dr' <0 $$
 for $0 <r \leq r_{\infty}$ and
 $$\frac{du}{dr}=-\frac{\mathsf{f}(1)}{3}r+O(r^2)$$
 as $r \rightarrow +0$. Therefore, considering this situation, we put
 
 \begin{Definition}
 Let $0<r_{\infty}$ and let $u\in\mathfrak{E}$ be an admissible solution with parameters
 $\mathsf{b}, u_{\mathsf{O}}$. If 
 \begin{equation}
\frac{\partial u}{\partial r} \leq -\frac{1}{C}r 
\quad\mbox{for}\quad 0<r\leq r_{\infty}, |\zeta|\leq 1\label{mntn}
\end{equation}
with a positive constant $C$, then the admissible solution $u$ said to be monotone.
 \end{Definition}

\begin{Remark}
Of course, if \eqref{mntn} holds with a positive constant $C$, the condition {\bf (a1)} holds for arbitrarily small positive $r_0$.
\end{Remark}

We claim
\begin{Theorem}\label{Th.3}
Let $0<r_{\infty}$. Let $\bar{u}\in \mathfrak{E}$ be a monotone
admissible solution with parameters $\mathsf{b}, \bar{u}_{\mathsf{O}}$ which satisfies the condition {\bf (HL)}. Then the solution $u$  of Theorem \ref{Th.1}
are monotone, provided that $\delta, \epsilon$ are sufficiently small, and 
\begin{equation}
|\nabla(\mathfrak{g}^{\flat}-\bar{\mathfrak{g}}^{\flat})|\leq  \delta'\cdot r,\label{gg}
\end{equation}
$\delta'$ being a sufficiently small positive number.
\end{Theorem}

\begin{Remark}
Note that 
\begin{align}
\nabla \mathfrak{g}^{\flat}&=
\frac{d\mathsf{b}}{d\varpi}\Big(
\frac{x_1}{\varpi}\frac{\partial}{\partial x_1}+
\frac{x_2}{\varpi}\frac{\partial}{\partial x_2}\Big) \nonumber \\
&=u_{\mathsf{O}}^{-1}\Omega(\mathsf{a}\varpi)^2\varpi
\Big(
\frac{x_1}{\varpi}\frac{\partial}{\partial x_1}+
\frac{x_2}{\varpi}\frac{\partial}{\partial x_2}\Big).
\end{align}

Thus we have
\begin{equation}
|\nabla(\mathfrak{g}^{\flat}-
\bar{\mathfrak{g}}^{\flat})|
\leq \|u_{\mathsf{O}}^{-1}\Omega^2-
\bar{u}_{\mathsf{O}}^{-1}\bar{\Omega}^2\|_{L^{\infty}}\cdot r,
\end{equation} 
and we see that \eqref{gg} holds with a small $\delta'>0$ if 
$ \|u_{\mathsf{O}}^{-1}\Omega^2-
\bar{u}_{\mathsf{O}}^{-1}\bar{\Omega}^2\|_{L^{\infty}} $ 
is sufficiently small.
\end{Remark}

In order to prove Theorem \ref{Th.3}, we use the following

\begin{Lemma}\label{Lem.3}
If $f\in\mathfrak{E}$ satisfies
$f^{\flat}\in C_0^1(B(r_{\infty}))$, then
\begin{equation}
|\nabla(\mathcal{K}f)^{\flat}(x)|\leq C\|f\|_1\cdot r,
\end{equation}
where
$$\|f\|_1=\|f\|_{\mathfrak{E}}+\|\nabla f^{\flat}\|_{L^{\infty}}.$$
\end{Lemma}

Proof. It is well known that 
$$\|\mathcal{K}f\|_{\mathfrak{E}}\leq C\|f\|_{\mathfrak{E}}. $$
Moreover, since
$$\frac{\partial}{\partial x_i}(\mathcal{K}f)^{\flat}({x})=
-\frac{1}{4\pi}\int
\frac{x_i-(x_i)'}{|{x}-{x}'|^3}
f^{\flat}({x}')d{x}',$$
we have
$$\Big|\frac{\partial}{\partial x_i}(\mathcal{K}f)^{\flat}\Big|\leq C\|f\|_{\mathfrak{E}}.$$
But an integration by parts gives
$$\partial_{i}((\mathcal{K}f)^{\flat})=-\mathcal{K}^{\flat}(\partial_{i}f^{\flat}),$$
where $\partial_i$ stands for $\displaystyle\frac{\partial}{\partial x_i}$,
so that
$$\Big|\frac{\partial^2}{\partial x_j\partial x_i}(\mathcal{K}f)^{\flat}({x})\Big|
=|-\partial_{j}\mathcal{K}^{\flat}(\partial_{i}f)^{\flat}|\leq C\|f\|_{1}.$$
On the other hand, 
$$\partial_{i}(\mathcal{K}f)^{\flat}({x})=-
\partial_{i}(\mathcal{K}f)^{\flat}(-{x}) $$
implies
$$\partial_{i}(\mathcal{K}f)^{\flat}(O)=0.$$
Therefore we get
$$|(\partial_{i}(\mathcal{K}f)^{\flat}({x})|\leq C\|f\|_{1}|x^j|.$$
This completes the proof of Lemma. $\square$\\


Proof of Theorem \ref{Th.3}: Look at
\begin{align*}
&u^{\flat}=\mathfrak{g}^{\flat}+1+(\mathcal{K}(\mathsf{f}\circ u))^{\flat}
-(\mathcal{K}(\mathsf{f}\circ u))^{\flat}(O), \\
&\bar{u}^{\flat}=\bar{\mathfrak{g}}^{\flat}+1+
(\mathcal{K}(\mathsf{f}\circ \bar{u}))^{\flat}
-(\mathcal{K}(\mathsf{f}\circ \bar{u}))^{\flat}(O).
\end{align*}
Therefore, with arbitrarily small $\epsilon'$, 
$$|\nabla(u^{\flat}-\bar{u}^{\flat})|\leq C\epsilon'\cdot r,
$$
since Lemma \ref{Lem.3} can be applied for $$\|\mathsf{f}\circ u^{\flat}-\mathsf{f}\circ \bar{u}^{\flat}\|_1\leq C'\epsilon' $$
and the condition \eqref{gg} is supposed. This implies
$$\Big|\frac{\partial u}{\partial r}-\frac{\partial\bar{u}}{\partial r}\Big|\leq \epsilon''\cdot r $$
so that
$$\frac{\partial u }{\partial r}\leq -\Big(
\frac{1}{\bar{C}}-\epsilon''\Big)\cdot r,$$
while we are supposing 
$$\frac{\partial\bar{u}}{\partial r}\leq -\frac{1}{\bar{C}}r\quad
\mbox{for}\quad 0<r\leq r_{\infty}, |\zeta|\leq 1.$$

Taking $\epsilon''\leq 1/2\bar{C}$, we have
$$\frac{\partial u}{\partial r}\leq -\frac{1}{2\bar{C}}\cdot r.
$$
This completes the proof. $\square$

\begin{Corollary}\label{Cor.2}
Suppose the same assumption as Corollary \ref{Cor.1}. Then we can claim that the solution $u$ is monotone, provided that
\begin{equation}
|\nabla \mathfrak{g}^{\flat}|\leq \delta'\cdot r,
\end{equation}
$\delta'>0$ being sufficiently small, or in other words, provided that 
$A_1u_{\mathsf{O}}^{-\nu}\|\Omega^2\|_{L^{\infty}(0,r_{\infty})}$ is sufficiently small.
\end{Corollary}

\section{Oblateness of the surface}

Let us assume that the equation of state is the exact $\gamma$-law:
$P=\mathsf{A}\rho^{\gamma}, 6/5<\gamma <2$, and let us denote
the Lane-Emden function of index $\nu$ by
$\theta(r;\nu)$. To fix the idea we put
$$\theta(r;\nu)=-\mu_1(\nu)
\Big(\frac{1}{\xi_1(\nu)}-\frac{1}{r}\Big)\quad
\mbox{for}\quad r \geq \xi_1(\nu),
$$
while 
$$\theta(r;\nu)>0\quad\mbox{for}\quad 0\leq r<\xi_1(\nu).$$
Here
$$\mu_1(\nu):=\int_0^{\xi_1(\nu)}
\theta(r,\nu)^{\nu}r^2dr
=-r^2\frac{d\theta(r;\nu)}{dr}\Big|_{r=\xi_1(\nu)}>0.$$

Thanks to the Corollary \ref{Cor.1}, we can consider slowly rotating configurations near to $\bar{u}=\theta(r;\nu), 1<\nu <5$. Let us assume the angular velocity $\Omega$ is a constant, and put
\begin{equation}
\beta=\frac{\Omega^2}{2\pi\mathsf{G}}
\Big(\frac{\mathsf{A}\gamma}{\gamma-1}\Big)^{\nu}u_{\mathsf{O}}^{-\nu}.
\end{equation}
Then 
\begin{equation}
\mathsf{b}(\varpi)=\frac{\beta}{4}\varpi^2,
\qquad
\mathfrak{g}(r,\zeta)=\frac{\beta}{4}r^2(1-\zeta^2).
\end{equation}
In our previous work \cite{JJTM} we constructed the distorted Lane-Emden function
$u=\Theta(r,\zeta;\nu,\beta)$ for small $\beta$, provided that $2\leq \nu <5$. But thanks to Corollary \ref{Cor.1} of the present study the distorted Lane-Emden function is defined even for $1<\nu <2$. Let us consider these distorted Lane-Emden configurations.\\

In \cite{JJTM} we showed the oblateness of the surface of the distorted configuration:
we proved that 
\begin{equation}
\sigma:=\frac{\Xi_1(0;\nu,\beta)-\Xi_1(\pm 1;\nu,\beta)}{\xi_1(\nu)}
\end{equation}
turns out to be positive for small $\beta$, provided that $2\leq \nu <5$. 
Here
$r=\Xi_1(\zeta; \nu,\beta)$ is the curve of the vacuum boundary, that is,
for $0\leq r\leq r_{\infty}$, $r_{\infty}$ being fixed so that 
$r_{\infty}>\xi_1(\nu)$, it holds
$$\Theta(r,\zeta;\nu,\beta)>0\quad
\Leftrightarrow \quad 0\leq r<\Xi_1(\zeta;\nu,\beta).
$$

We are going to prove that it is the case even if $1<\nu <2$.\\

Suppose $1<\nu <2$. Let us put
\begin{equation}
\mathfrak{g}(r,\zeta)=\beta\mathfrak{g}_1(r,\zeta),\quad
\mathfrak{g}_1(r,\zeta)=\frac{1}{4}r^2(1-\zeta^2).
\end{equation}
The distorted Lane-Emden function $\Theta$ is the solution of
\begin{equation}
\Theta=\beta\mathfrak{g}_1+\mathcal{G}(\Theta).
\end{equation}
On the other hand, Lemma \ref{Lem.1} implies the following

\begin{Proposition}
Let $1<\nu<\infty$. If $u,u+h \in \mathfrak{E}$ satisfy
$\|u\|_{\mathfrak{E}}, \|u+h\|_{\mathfrak{E}}\leq M$, then
\begin{equation}
\|\mathcal{G}(u+h)-\mathcal{G}(u)
-D\mathcal{G}(u)h\|_{\mathfrak{E}}
\leq C\|h|_{\mathfrak{E}}^{\nu\wedge 2}.
\end{equation}
Here the constant $C$ depends upon $M$.
\end{Proposition}
Therefore we can claim
\begin{equation}
\Theta-\theta=\beta\mathfrak{g}_1+
D\mathcal{G}(\theta)(\Theta-\theta)+
O_{\mathfrak{E}}(\|\Theta-\theta\|_{\mathfrak{E}}^{\nu}),
\end{equation}
where the symbol $O_{\mathfrak{E}}(\kappa)$ stands for functions $\varphi
\in \mathfrak{E}$ such that $\|\varphi\|_{\mathfrak{E}}\leq C\kappa$.
Since $(I-D\mathcal{G}(\theta))^{-1}$ exists to be a bounded linear operator in $\mathfrak{E}$, we have
\begin{equation}
\|\Theta-\theta\|_{\mathfrak{E}}\leq C\|\beta\mathfrak{g}_1\|_{\mathfrak{E}}\leq C'\beta,
\end{equation}
keeping in mind that $\nu >1$, so 
\begin{equation}
\Theta=\theta+\beta\mathfrak{h}+O_{\mathfrak{E}}(\beta^{\nu}),
\end{equation}
where
\begin{equation}
\mathfrak{h}:=(I-D\mathcal{G}(\theta))^{-1}\mathfrak{g}_1.
\end{equation}
Thus we are going to observe the structure of the function
$\mathfrak{h}$. 

Put
\begin{equation}
h_j(r)=\frac{2j+1}{2}\int_{-1}^1\mathfrak{h}(r,\zeta)P_j(\zeta)d\zeta.
\end{equation}
Here $P_j$ is the Legendre's polynomial. Since $\mathfrak{h}$ satisfies
\begin{equation}
\mathfrak{h}=D\mathcal{G}(\theta)\mathfrak{h}+\mathfrak{g}_1,
\end{equation}
we can assume that $\mathfrak{h}$ is defined for $0\leq r<+\infty, |\zeta|\leq 1$ so that
\begin{equation}
\mathfrak{h}-\mathfrak{g}_1=O(\frac{1}{r}),\qquad
\frac{\partial}{\partial r}(\mathfrak{h}-\mathfrak{g}_1)=O(\frac{1}{r^2})
\end{equation}
as $r\rightarrow +\infty$, since the support of 
$\mathsf{f}'(\theta)^{\flat}=[\nu(\theta\vee 0)^{\nu-1}]^{\flat}$ is $B(\xi_1(\nu))$. Note that
\begin{equation}
\mathfrak{g}_1(r,\zeta)=\frac{r^2}{6}-
\frac{r^2}{6}P_2(\zeta),
\qquad
P_2(\zeta)=\frac{1}{2}(3\zeta^2-1).
\end{equation}\\

Since $\mathfrak{h}$ is equatorially symmetric, $h_j$ vanishes for any odd $j$. \\

Let $j$ be an even number such that $j\geq 4$. Then $h_j$ is the solution of
the equation
$$\Big[-\frac{1}{r^2}\frac{d}{dr}r^2\frac{d}{dr}+
\frac{j(j+1)}{r^2}\Big]h_j=\nu(\theta\vee 0)^{\nu-1}h_j $$
such that
$h_j=O(1/r), dh_j/dr =O(1/r^2)$ as $r\rightarrow +\infty$. Therefore we can claim that $h_j$ vanishes by the argument in the proof of 
Theorem \ref{Th.2}.\\

Let us consider $j=2$. Then $h_2$ is the solution of
the equation
$$\Big[-\frac{1}{r^2}\frac{d}{dr}r^2\frac{d}{dr}+
\frac{6}{r^2}\Big]h_2=\nu(\theta\vee 0)^{\nu-1}h_2 $$
such that
\begin{equation}
h_2(r)=-\frac{r^2}{6}+O\Big(\frac{1}{r}\Big),
\qquad
\frac{d}{dr}h_2(r)=-\frac{r}{3}
+O\Big(\frac{1}{r^2}\Big)
\end{equation}
as $r\rightarrow +\infty$. Then we have the integral representation
\begin{equation}
h_2(r)=\frac{r^2}{6}
\Big[-1+\frac{6}{5}\frac{1}{r^5}
\int_0^r\mathsf{f}'(\theta(r'))h_2(r')r'^4dr'+
\frac{6}{5}r^3\int_r^{\infty}\mathsf{f}'(\theta(r'))h_2(r')(r')^{-1}dr'\Big] \label{5.16}
\end{equation}
where $\mathsf{f}'(\theta(r))=\nu(\theta(r;\nu)\vee 0)^{\nu-1}$. 
 See Appendix, namely, note that we apply Proposition \ref{Prop.A} to $j=2, A=-5/6$.  We are going to prove 
the following

\begin{Theorem}\label{Th.4}
Let $1<\nu <2$. It holds that
\begin{equation}
h_2(r)<0\quad\mbox{for}\quad 0<r\leq\xi_1(\nu).
\end{equation}
\end{Theorem}

Proof. Since the equation
$$\Big[-\frac{1}{r^2}\frac{d}{dr}r^2\frac{d}{dr}+\frac{6}{r^2}\Big]y=
\nu(\theta\vee 0)^{\nu-1}y
$$
admits a fundamental system of solutions
$y_1 \sim r^2, y_2 \sim r^{-3}$ as $r\rightarrow +0$, we see
$h_2(r)=O(r^2)$ as $r\rightarrow +0$. Then the integral representation \eqref{5.16}
implies
\begin{equation}
h_2(r)=-\frac{r^2}{6}(1+O(r^2))\quad\mbox{as}\quad r\rightarrow +0.\label{5.18}
\end{equation}

Let us consider
\begin{equation}
H(r):=\frac{h_2(r)}{\psi(r)},
\end{equation}
where
\begin{equation}
\psi(r):=\frac{1}{r^2}\int_0^r
\theta(r')^{\nu}r'^2dr'=
-\frac{d\theta}{dr} >0.
\end{equation}
Then \eqref{5.16} reads
\begin{align}
H(r)&=\frac{1}{5\psi(r)}\Big[ -\frac{5}{6}r^2+
\frac{1}{r^2}\int_0^r(-D\theta(r'))H(r')\Big(\frac{r'}{r}\Big)r'^3dr'+ \nonumber \\
&+r\int_r^{\xi_1}(-D\theta(r'))H(r')\Big(\frac{r}{r'}\Big)dr'\Big] \label{5.19}
\end{align}
for $0<r\leq \xi_1=\xi_1(\nu)$, where
$$D\theta(r)=\frac{d\theta(r;\nu)}{dr}.$$

Since 
$$\psi(r)=\frac{r}{3}(1+O(r^2))\quad\mbox{as}\quad r\rightarrow +0,
$$
\eqref{5.18} implies
$$H(r)=-\frac{r}{2}(1+O(r^2))<0\quad \mbox{as}\quad r\rightarrow +0.
$$

We want to show 
$$H(r)<0\quad\mbox{for}\quad 0<r\leq \xi_1.
$$
Otherwise, there would exist $r_*\in ]0,\xi_1]$ such that
$$H(r_*)=\max_{0\leq r\leq \xi_1}H(r) \geq 0.$$
Then \eqref{5.19} would imply
\begin{align*}
H(r)&\leq \frac{1}{5\psi(r)}\Big[-\frac{5}{6}r^2+ \\
&+\Big(\frac{1}{r^2}\int_0^r(-D\theta(r'))r'^3dr'+
r\int_r^{\xi_1}(-D\theta(r'))dr'\Big)H(r_*)\Big] \\
&=\frac{1}{5\psi(r)}\Big[-\frac{5}{6}r^2+
\frac{3}{r^2}\psi(r)H(r_*)\Big] \\
&=-\frac{r^2}{6\psi(r)}+\frac{3}{5}H(r_*)
\end{align*}
for $0<r\leq \xi_1$. Thus
$$H(r_*)\leq -\frac{r_*^2}{6\psi(r_*)}+\frac{3}{5}H(r_*),
$$
so that
$$H(r*)\leq -\frac{5r_*^2}{12\psi(r_*)}<0,$$
a contradiction. This completes the proof. $\square$\\

Summing up, we have
\begin{equation}
\mathfrak{h}(r,\zeta)=h_0(r)+h_2(r)P_2(\zeta)
\end{equation}
and
\begin{equation}
h_2(r)<0\quad\mbox{for}\quad 0<r\leq \xi_1(\nu).
\end{equation}

Then it follows that
\begin{equation}
\Xi_1(\zeta;\nu,\beta)=
\xi_1(\nu)+
\frac{\xi_1(\nu)^2}{\mu_1(\nu)}\mathfrak{h}(\xi_1(\nu),\zeta)\beta+
O(\beta^{\nu})
\end{equation}
and
\begin{equation}
\sigma=-\frac{3}{2}\frac{\xi_1(\nu)}{\mu_1(\nu)}
h_2(\xi_1(\nu))\beta+
O(\beta^{\nu}) >0,
\end{equation}
provided that $0<\beta \ll 1$,
since $\xi_1(\nu)>0, \mu_1(\nu)>0$
and $-h_2(\xi_1(\nu))>0$. This was to be demonstrated.

\section{Relation between the total mass and the central density}

In this section we discuss about the relation between the total mass and the central density of the models of uniformly rotating stars constructed for the exact $\gamma$-law: $P=\mathsf{A}\rho^{\gamma}$ with $ 1<\nu=\frac{1}{\gamma-1} <5$.\\

The density distribution is given by
\begin{equation}
\rho=\rho_{\mathsf{O}}(\Theta(r/\mathsf{a},\zeta;\nu,\beta)\vee 0)^{\nu}, \quad
\mathsf{a}=\sqrt{\frac{\mathsf{A}\gamma}{4\pi\mathsf{G}(\gamma-1)}}
\rho_{\mathsf{O}}^{-\frac{2-\gamma}{2}}, \quad
\beta=\frac{\Omega^2}{2\pi\mathsf{G}\rho_{\mathsf{O}}},
\end{equation}
provided that $\beta$ is sufficiently small.

The total mass
\begin{equation}
 M: =2\pi\int_{-1}^1\int_0^{\infty} \rho(r,\zeta)r^2drd\zeta 
\end{equation}
is given by
\begin{equation}
M=\rho_{\mathsf{O}}\mathsf{a}^3M_1(\nu,\beta)=
\Big(\frac{\mathsf{A}\gamma}{4\pi\mathsf{G}(\gamma-1)}\Big)^{3/2}
\rho_{\mathsf{O}}^{\frac{-4+3\gamma}{2}}
M_1(\nu,\beta)
\end{equation}
with
\begin{equation}
M_1(\nu,\beta):=2\pi\int_{-1}^1\int_0^{\infty}
(\Theta(r,\zeta;\nu,\beta)\vee 0)^{\nu}r^2drd\zeta.
\end{equation}

Since
\begin{align*}
\Big(\frac{\partial M}{\partial\rho_{\mathsf{O}}}\Big)_{\Omega=\mbox{Const.}}&=
\Big(\frac{\mathsf{A}\gamma}{4\pi\mathsf{G}(\gamma-1)}\Big)^{3/2}
\rho_{\mathsf{O}}^{\frac{-4+3\gamma}{2}-1}
\Big(\frac{-4+3\gamma}{2}M_1-\beta\frac{\partial M_1}{\partial\beta}\Big) \\
&\not=0
\end{align*}
if $-4+3\gamma \not=0$ and $\beta$ is sufficiently small, 
we have the following

\begin{Theorem}
 Assume the exact $\gamma$-law, $P=\mathsf{A}\rho^{\gamma}$, with 
$\frac{6}{5}<\gamma <2$, and consider small constant angular velocities $\Omega$. Then
the central density $\rho_{\mathsf{O}}$ is uniquely determined by the total mass $M$ for each
given angular velocity $\Omega$, provided that $\gamma \not=4/3$ and $\beta$ is
sufficiently small.
\end{Theorem}

\begin{Remark}
When $\gamma=4/3$, the situation is too delicate. We should investigate the behavior of $\partial M_1/\partial\beta$ thoroughly.
\end{Remark}

Note that the above conclusion can be expressed as follows:

Let us consider the exact $\gamma$-law $P=\mathsf{A}\rho^{\gamma}$ with $\frac{6}{5}<\gamma <2, \gamma \not=\frac{4}{3}$. The total mass $M=\mathcal{M}(
\rho_{\mathsf{O}}, \Omega^2)$ is determined for the unique configuration of the central density $\rho_{\mathsf{O}}$ and the angular velocity 
$\Omega$. When we fix $\bar{M}=\mathcal{M}(\bar{\rho}_{\mathsf{O}}, 0)$, then we can find a continuous curve $\mathcal{C}: \Omega^2 \mapsto \rho_{\mathsf{O}}$ 
defined for $\Omega^2 \ll 1$ such that $\mathcal{C}(0)=\bar{\rho}_{\mathsf{O}}$ and $$\mathcal{M}(\mathcal{C}(\Omega^2), \Omega^2)=\bar{M}.$$\\

It is easy to prove by the implicit function theorem that the same situation holds for any general equation of state $\rho \mapsto P$, which is not necessarily the exact 
$\gamma$-law, as follows, thanks to Theorem 2:

\

{\bf Theorem 5bis}
{\it Suppose that the assumption of Theorem 2 holds, and consider the central density 
$\rho_{\mathsf{O}}$ such that
$|\rho_{\mathsf{O}}-\bar{\rho}_{\mathsf{O}}| \ll 1$, where
$$\bar{\rho}_{\mathsf{O}}=
\Big(\frac{\gamma-1}{\mathsf{A}\gamma}\Big)^{\nu}\bar{u}_{\mathsf{O}}^{\nu}(1+\Lambda_{\rho}(\bar{u}_{\mathsf{O}})),$$
and small constant angular velocities $\Omega$. Let $\mathcal{M}(\rho_{\mathsf{O}}, \Omega^2)$ be the total mass of the unique configuration of the central density $\rho_{\mathsf{O}}$ and the angular velocity $\Omega$. Then if we assume
\begin{equation}
\frac{\partial}{\partial\rho_{\mathsf{O}}}\mathcal{M}({\rho}_{\mathsf{O}}, 0)\Big|_{\rho_{\mathsf{O}}=\bar{\rho}_{\mathsf{O}}}\not=0,
\end{equation}
then we have a curve $\mathcal{C}: \Omega^2 \mapsto \rho_{\mathsf{O}}$ such that
$\mathcal{C}(0)=\bar{\rho}_{\mathsf{O}}$ and $$\mathcal{M}(\mathcal{C}(\Omega^2),\Omega^2)=
\mathcal{M}(\bar{\rho}_{\mathsf{O}},0).$$
}

(Cf. \cite[Theorem 2.2]{StraussW}. )


\section{Supplementary Remark}


In this article we have constructed solutions under the assumption that the angular velocity $\Omega$ is a prescribed function of the distance $\varpi$ from the rotating axis. However in many astrophysical literatures it is assumed that $j=\Omega \varpi^2$ is a prescribed function  not directly of  $\varpi$ but of the quantity
$$m=\int_{\varpi(x)\leq \varpi}\rho(x)dx.$$
For example, J. P. Ostriker and J. W-K. Mark in 1968 wrote
\begin{quote}
The centrifugal force is found from the equilibrium distribution of angular momentum; the latter, unfortunately, is indeterminate. However, Mistel (1963) --\cite{Mistel}-- and Crampin and Hyle
(1964) --\cite{CrampinH}-- have pointed out that the angular momentum of a given toroidal element will be conserved throughout the evolution of a system; i. e., $j=j(m)$ is an
invariant function, if $m$ is a Lagrangian coordinate denoting a mass element.
(\cite[p. 1077]{OstrikerM})
\end{quote}

Therefore let us consider the case in which $j$ is a given function in $C^1([0,+\infty))$ such that
$j(0)=0$ 
and the quantities $\Omega$ and so on
are defined as follows:

\begin{align}
\Omega(\varpi)&=j(m(\varpi))\varpi^{-2}, \\
m(\varpi)&=2\pi\int_{-\infty}^{+\infty}\int_0^{\varpi}
\rho(r',\zeta')\varpi'd\varpi'dz' \nonumber \\
&=\int_{\varpi(x)\leq \varpi}\rho^{\flat}(x)dx, \\
&\mbox{with} \quad r'=\sqrt{(\varpi')^2+(z')^2},\quad \zeta'=\frac{z'}{r'}, \nonumber \\
\mathsf{B}(\varpi)&=\int_0^{\varpi}\Omega(\varpi')^2\varpi'd\varpi' \nonumber \\
&=\int_0^{\varpi}j(m(\varpi'))^2(\varpi')^{-3}d\varpi', \label{SR3}\\
\mathsf{b}(\varpi)&=u_{\mathsf{O}}^{-1}\mathsf{B}(\mathsf{a}\varpi), \label{SR4}\\
\mathfrak{g}(r,\zeta)&=\mathsf{b}(r\sqrt{1-\zeta^2}).\label{SR5}
\end{align}

\begin{Remark}

The quantity $j=j(m)=\Omega(\varpi)\varpi^2$ can be called
{\bf `angular momentum per unit mass'}, as in e.g., \cite{OstrikerM}, in the following sense.

The angular momentum for an Eulerian  gaseous flow is defined as
$$\vec{J}=\int \vec{x}\mathsf{x}(\rho\vec{v})d\vec{x},$$
where $\vec{x}=(x_1,x_2,x_3)$ is the static Cartesian coordinates. When $\rho=\tilde{\rho}
(\varpi, z)$, and
$$\vec{v}=-\Omega(\varpi)x_2\frac{\partial}{\partial x_1}+
\Omega(\varpi)x_1\frac{\partial}{\partial x_2}, $$
then we have $\vec{J}=(0,0,J)^T$ with
$$J=2\pi\int_{-\infty}^{+\infty}\int_0^{+\infty}\tilde{\rho}(\varpi,z)
\Omega(\varpi)\varpi^3d\varpi dz.$$
On the other hand, since
$$m(\varpi)=2\pi\int_{-\infty}^{+\infty}\int_0^{\varpi}
\tilde{\rho}(\varpi',z)\varpi'd\varpi'dz,$$
we have 
$$dm=2\pi\tilde{\rho}(\varpi, z)\varpi d\varpi dz,
$$
so that
$$J=\int_0^M\Omega \varpi^2 dm=\int_0^Mj(m)dm,$$
where 
$M=m(+\infty) $ is the total mass, and $\varpi=\varpi(m)$ is a monotone increasing function of $m \in [0,M]$.
\end{Remark}

Then for $u \in \mathfrak{E}$, $r_{\infty}$ being fixed, the quantity $m(\varpi)=
(\mathcal{M}u)(\varpi)$ is given by
\begin{equation}
m(\varpi)=
2\pi(4\pi\mathsf{G})^{-\frac{3}{2}}
\Big(\frac{\mathsf{A}\gamma}{\gamma-1}\Big)^{\frac{1}{2(\gamma-1)}}
u_{\mathsf{O}}^{\frac{3\gamma-4}{2(\gamma-1)}}\int_{-\infty}^{+\infty}
\int_0^{\varpi/\mathsf{a}}\mathsf{f}(u(r',\zeta'))\varpi'd\varpi'dz', \label{SR6}
\end{equation}
where $r'=\sqrt{(\varpi')^2+(z')^2},\quad \zeta'=\frac{z'}{r'}$. Of course the function
 $\mathsf{f}$ is defined by \eqref{Def.f}.

As the functional space for $j$, let us denote by $\mathfrak{J}$ the Banach space of functions $j \in C^1[0,+\infty)$ such that $j(0)=0$ and
\begin{equation}
\|j\|_{\mathfrak{J}}:=
\sup_{m\in [0,+\infty)}|j(m)|+
\sup_{m\in [0,\infty)}\Big|\frac{dj}{dm}\Big| <\infty.
\end{equation}

Let us define the operator $\mathcal{B}$ by
\begin{equation}
\mathcal{B}(j, u)(r,\zeta)=\mathfrak{g}(r,\zeta)=\mathsf{b}(r\sqrt{1-\zeta^2})
\end{equation}
in view of \eqref{SR6},\eqref{SR3},\eqref{SR4},\eqref{SR5}.
This $\mathcal{B}$ is a nonlinear mapping from $\mathfrak{J}\times \mathfrak{E} $ into $\mathfrak{E}$ and the 
Fr\'{e}chet partial derivative $D_u\mathcal{B}(j,u)$ with respect to $u$ is given by
\begin{equation}
(D_u\mathcal{B}(j, u).h)(r,\zeta)=
u_{\mathsf{O}}^{-1}
\int_0^{\mathsf{a}r\sqrt{1-\zeta^2}}
2j\frac{dj}{dm}\Big|_{m=\mathcal{M}(u)(\varpi')}\cdot (D\mathcal{M}(u).h)(\varpi')
(\varpi')^{-3}d\varpi',
\end{equation}
where
\begin{align}
(D\mathcal{M}(u).h)(\varpi)&=
2\pi(4\pi\mathsf{G})^{-\frac{3}{2}}
\Big(\frac{\mathsf{A}\gamma}{\gamma-1}\Big)^{\frac{1}{2(\gamma-1)}}
u_{\mathsf{O}}^{\frac{3\gamma-4}{2(\gamma-1)}}\times \nonumber \\
&\times
\int_{-\infty}^{+\infty}\int_0^{\varpi/\mathsf{a}}
\mathsf{f}'(u(r',\zeta'))h(r',\zeta')\varpi'd\varpi'dz'.
\end{align}
Of course $\mathsf{f}'(u)$ is that given by
\eqref{2.25},\eqref{2.26}. 
By Ascoli-Arcela's theorem we see that $D_u\mathcal{B}(j,u)$ is a compact operator in
$\mathfrak{E}$. The mapping $(j,u) \mapsto D_u\mathcal{B}(j,u)$ is continuous 
re the norms of $\mathfrak{J}, \mathfrak{E}$ and the operator norm of bounded linear operators 
in $\mathfrak{E}$. 
Let us introduce the convenient notation $$\mathsf{b}(\varpi)=\mathcal{B}^{\sharp}(j,u)(\varpi)=\mathcal{B}^{\sharp}(j,u)(r\sqrt{1-\zeta^2})=
\mathcal{B}(j,u)(r,\zeta).$$

In this framework, Theorem 1 holds if we replace the condition {\bf (HL)} by the following {\bf (mHL)}:
\begin{equation}
h\in\mathfrak{E}\quad
h=[D_u\mathcal{B}(j, u)+D\mathcal{G}(u)]h\quad
\Rightarrow \quad h=0.
\end{equation}

Precisely speaking, we can claim the following \\

{\bf Theorem 1S } {\it Let $r_{\infty}>0$. Let $\bar{u}\in\mathfrak{E}$ be an admissible solution of \eqref{E} with parameters
$\bar{\mathsf{b}}=\mathcal{B}^{\sharp}(\bar{j},\bar{u}), \bar{u}_{\mathsf{O}}$. Here $\bar{j}\in\mathfrak{J}$ and $u_{\mathsf{O}}>0$. Suppose that $u=\bar{u}, j=\bar{j}$ enjoys the condition {\bf (mHL)}.
Then there exist small positive constants $\delta,\epsilon$ such that
for $(j,u_{\mathsf{O}}) \in \mathfrak{J}\times \mathbb{R}$ such that
$\|j-\bar{j}\|_{\mathfrak{J}}+|u_{\mathsf{O}}-\bar{u}_{\mathsf{O}}|\leq
\delta$ there is a unique admissible solution $u$ with parameters
$\mathsf{b} =\mathcal{B^{\sharp}}(j,u), u_{\mathsf{O}}$ which enjoys  {\bf (mHL)} such that
$\|u-\bar{u}\|_{\mathfrak{E}}\leq\epsilon$ and $(j,u_{\mathsf{O}}) \mapsto u$ is continuous. }


\vspace{10mm}

{\bf\large Acknowledgment}\\

The first author acknowledges the support from NSF Grant  DMS-1608494. The second author  expresses his sincere thanks to Professor 
Todd A. Oliynik ( Monash University) who kindly drew his attention to the work by U. Heilig, and to the support from Yamaguchi University Foundation Research Grant A1-2 (2017)  and the support from
JSPS KAKENHI Grant JP 18K0337. 

\vspace{10mm}

{\bf\large Appendix}\\

Here we give an elementary proof of the integral representation 
which we have used in the discussion on the condition {\bf (HL)} and oblateness.

\begin{Proposition}\label{Prop.A}
 Let $j$ be a positive integer, and $q$ be a continuous function on $[0,+\infty[$ such that $q(r)=0$ for $r\geq R$, $R$ being a finite positive number. Let $y=y(r)$ be a solution of the equation
$$\Big[-\frac{1}{r^2}\frac{d}{dr}
r^2\frac{d}{dr}+\frac{j(j+1)}{r^2}\Big]y=
q(r)y \qquad (r>0) \leqno (E)$$
such that 
$$y=O(1)\quad\mbox{as}\quad r\rightarrow +0 \leqno(B_0)$$
and
$$r^{-j}\Big(r\frac{dy}{dr}+(j+1)y\Big) \rightarrow A
\quad\mbox{as}\quad r\rightarrow +\infty. \leqno(B_{\infty})$$
Then there holds the integral representation
\begin{align*}
y(r)&=\frac{1}{2j+1}\Big[
A r^j+\frac{1}{r^2}\int_0^r
q(s)y(s)\Big(\frac{s}{r}\Big)^{j-1}s^3ds + \\
&+r\int_r^{+\infty}q(s)y(s)\Big(\frac{r}{s}\Big)^{j-1}ds\Big].
\end{align*}
\end{Proposition}

Proof. First we note that $(B_0)$ implies
$$y=O(r^j)\quad\mbox{as}\quad r\rightarrow +0, $$
since the ordinary differential equation $(E)$ has a fundamental system of solutions $y=\psi_1(r), \psi_2(r)$ such that
$\psi_1(r)\sim r^j, \psi_2(r)\sim r^{-j-1}$ as $ r\rightarrow +0$.

Now the left-hand side of $(E)$ can be written as
$$-r^{j-1}\frac{d}{dr}\Big(r^{-2j}\frac{d}{dr}(r^{j+1}y)\Big). $$
Therefore the integration of the equation gives
$$r^{-2j}\frac{d}{dr}(r^{j+1}y)=
\int_r^{+\infty}q(s)y(s)s^{-j+1}ds + C,
$$
$C$ being a constant. But, since we are supposing $(B_{\infty})$:
$$r^{-2j}\frac{d}{dr}(r^{j+1}y)=
r^{-j+1}\frac{dy}{dr}+(j+1)r^{-j}y \rightarrow A \quad\mbox{as}\quad r\rightarrow +\infty,$$ and since
$$\int_r^{+\infty}q(s)y(s)s^{-j+1}ds =0$$
for $r\geq R$, we see $C=A$, that is,
$$r^{-2j}\frac{d}{dr}(r^{j+1}y)=
\int_r^{+\infty}q(s)y(s)s^{-j+1}ds + A,
$$
that is,
$$\frac{d}{dr}(r^{j+1}y)=
r^{2j}\int_r^{+\infty}q(s)y(s)s^{-j+1}ds + Ar^{2j}.
$$
Integrating this equality, keeping in mind that $r^{j+1}y=O(r^{2j+1})
\rightarrow 0$ as $r\rightarrow +0$, and $q(s)y(s)s^{-j+1}=O(s)$ as $s\rightarrow +0$, we can claim
$$r^{j+1}y=
\int_0^rt^{2j}
\Big(\int_t^{+\infty}q(s)y(s)s^{-j+1}ds\Big)dt +
\frac{A}{2j+1}r^{2j+1},$$
that is,
$$y=r^{-j-1}\int_0^rt^{2j}\Big(
\int_t^{+\infty}q(s)y(s)s^{-j+1}ds\Big)dt+
\frac{A}{2j+1}r^j.$$
Change of order of integration leads us to the claimed integral representation. $\square$

\end{document}